\documentclass[11pt,letterpaper]{amsart}
\usepackage{mathpazo} 
\usepackage{amsfonts}
\usepackage{amsmath}
\usepackage{amssymb}
\usepackage{mathrsfs}
\usepackage{mathtools}
\usepackage{amsthm}
\usepackage{bm}
\usepackage{bbm}

\usepackage{yhmath}
\usepackage{tikz-cd}
\usepackage[T1]{fontenc}
\usepackage[utf8]{inputenc}
\usepackage[english]{babel}
\usepackage{color}
\definecolor{darkred}{RGB}{203,65,84}
\definecolor{darkblue}{RGB}{70,130,180}
\definecolor{brown}{RGB}{139,69,19}

\newcommand{\mapping}[4] { \left\{
    \begin{array}{rcl}
      #1 &\rightarrow& #2\\
      #3 &\mapsto& #4
    \end{array}
  \right.  }
\usepackage[colorlinks=true, linkcolor=darkblue, citecolor=darkred, urlcolor=cyan]{hyperref}

\newcommand{\defeq}{\coloneqq}

\newcommand{\Aut}{\mathrm{Aut}}

\newcommand{\Hopf}{\mathrm{Hopf}}
\newcommand{\N}{\mathbb N}
\newcommand{\U}{\mathsf U}
\newcommand{\T}{\mathsf T}
\newcommand{\bfS}{\mathbf S}
\newcommand{\g}{\mathbf{g}}
\newcommand{\E}{\mathrm{E}}
\newcommand{\A}{\mathrm{A}}
\newcommand{\dvol}{\mathrm{d\,vol}}
\renewcommand{\d}{\mathrm{d}}
\newcommand{\hh}{{\mathbf H^2}}

\renewcommand{\leq}{\leqslant}

\renewcommand{\geq}{\geqslant}
\renewcommand{\epsilon}{\varepsilon}
\newcommand{\cH}{\mathcal H}
\newcommand{\ind}{{\scriptscriptstyle{\mathrm{ind}}}}
\title{Minimal surfaces with negative curvature in large dimensional spheres}

\author[M. Ancona]{Michele Ancona}
\address{M. Ancona: Universit\'e C\^ote d'Azur, CNRS,  LJAD,  France}
\email{michele.ancona@univ-cotedazur.fr}

\author[F. Labourie]{Fran\c cois Labourie}
\address{F. Labourie: EPF Lausanne, SB-SCI-FL,  CH-1015 Lausanne,  Switzerland}
\email{francois.labourie@epfl.ch}

\author[A. Roig-Sanchis]{Anna Roig-Sanchis}
\address{A. Roig-Sanchis: Universit\'e C\^ote d'Azur, CNRS,  LJAD,  France}
\email{anna.roig-sanchis@univ-cotedazur.fr}

\author[J. Toulisse]{J\'er\'emy Toulisse}
\address{J. Toulisse: Universit\'e C\^ote d'Azur, CNRS,  LJAD,  France}
\email{jeremy.toulisse@univ-cotedazur.fr}

\newtheorem*{maintheorem}{Main Theorem}
\newtheorem*{maincorollary}{Main Corollary}

\newtheorem{lemma}{Lemma}[section]

\newtheorem{theorem}[lemma]{Theorem}

\newtheorem{proposition}[lemma]{Proposition}
\theoremstyle{definition}
\newtheorem{definition}[lemma]{Definition}
\theoremstyle{remark}
\newtheorem{remark}[lemma]{Remark}

\thanks{F.~L.  and J.~T. acknowledge funding by the European Research Council under ERC-Advanced grant 101095722. J.~T. acknowledges the support of the Institut Universitaire de France. F.~L. is supported by the Swiss State Secretariat for Education, Research and Innovation (SERI): MB25.00031.}
\begin{document}

\maketitle

\begin{abstract}
In this paper,  we answer positively a question of S.-T. Yau by proving the existence of closed minimal surfaces with negative induced curvature in any sphere of large dimension. The proof follows the strategy of A. Song, applying it  to closed Riemann surfaces with large automorphism groups, and obtaining almost hyperbolic minimal surfaces.
\end{abstract}

\tableofcontents

\section{Introduction}

A classical result due to Robert Bryant \cite{bryant} states the nonexistence of closed minimal surfaces with constant negative induced curvature in round spheres. In \cite{yau}, Shing-Tung Yau asks the question of existence of closed minimal surfaces in spheres with negative induced curvature. In this paper, we give, as a corollary of our main result, a positive  answer to this question:

\begin{maincorollary}
There exists a closed minimal surface of negative curvature in any round sphere of large dimension.	
\end{maincorollary}

The proof follows the asymptotic strategy -- when the dimension goes to infinity  --  developed by Antoine Song in \cite{song}. In particular, we do not have any information on the dimension of the sphere for which it happens, although it is a well-known consequence of Gau\ss--Codazzi equations due to Blaine Lawson \cite[Proposition 1.5]{Lawson:1970aa}  that closed negatively curved surfaces cannot exist in $\mathbf{S}^3$.

The novelty lies in the fact that we refine Song's construction and  get the stronger result
 
\begin{maintheorem}
For each integer $n$ large enough there exists a closed  minimal surface $\Sigma_n$ in the round sphere of dimension $n$ such that for all integer $k$ 
$$
\lim_{n\to\infty}\left\Vert \kappa_n+8\right\Vert_{C^k}=0 
$$
where $\kappa_n$ is the curvature of $\Sigma_n$.

Furthermore, one can choose the sequence $(\Sigma_n)_{n\in \mathbb N}$ in such a way that, for each $n$, there exists a finite group $\Gamma_n$ in the isometry group of the $n$-sphere such that the quotient  $\Gamma_n\backslash\Sigma_n$ is a Riemann surface independent of $n$.
\end{maintheorem}

In other words, asymptotically, Bryant's result does not hold. Let us give an idea of the proof of our Main Theorem: in his groundbreaking work, Song developed a new strategy to study harmonic maps into large dimensional spheres that are equivariant under  unitary representations. Applying this general construction to the hyperbolic sphere with $3$ cusps, Song obtained in \cite[Theorem 0.3]{song} a sequence of minimal surfaces in spheres that  converges in the Benjamini--Schramm sense to the hyperbolic plane. The lack of smooth convergence in his theorem comes from the noncompactness of the cusped hyperbolic surface.

The main idea behind using a $3$-punctured sphere is that its Teichmüller space is reduced to a point. In this paper, we replace the use of the $3$-punctured sphere by that of an orbifold whose Teichmüller space is a point. Using the notion of induced representations, we adapt the approach of Song to the orbifold case.

\subsection*{Acknowledgements} We thank R\'emy Rodiac, Andrea Seppi, Peter Smillie  and Antoine Song for their interests and useful remarks. We also thank the referees for their helpful remarks.

\section{Harmonic maps, minimal surfaces and results of Song}

\subsection{Minimal surfaces}

Let $\Sigma$ be a closed surface and $(M,h)$ be a Riemannian manifold. The \emph{area} of a smooth map $f: \Sigma \to M$ is defined by
\[\A(f) = \int_\Sigma \dvol_{f^*h}~.\]
The map $f$ is called \emph{minimal} if it is a critical point of the area, namely if for any smooth variation $(f_t)_{t\in (-\epsilon,\epsilon)}$ we have
\[\left. \frac{\d}{\d t}\right\vert_{t=0} \A(f_t)=0~.\]

We recall a definition of Robert Gulliver \cite{Gulliver:1973}, see also Mario Micallef and Brian White for minimal maps \cite{Micallef:1995aa}. Let $P$ be a discrete subset of $\Sigma$, a \emph{branched  minimal immersion} with {\em branch locus} $P$ is a map $f: \Sigma \to M$  which is\begin{itemize}
	\item a minimal immersion away from  $P$,
	\item for every $p$ in $P$, there is a positive integer $d$, called \emph{the order of branching}, such that the map $f$ can be written in local coordinates centered at $p$ as $f(z) = (z^{d+1}, g(z))$, where  $g$ and $\mathrm{d}g$ are $\mathcal{O}(z^{d+1})$.
\end{itemize}
\begin{remark}\label{rem:SingularCurvature}
If $p$ is a branched point of a branched minimal immersion $f$, then the differential of $f$ vanishes at $p$. Hence, the induced metric, seen as a 2-tensor on $\Sigma$, vanishes at $p$.
\end{remark}

\subsection{Harmonic maps from surfaces}

Let $(\Sigma,g)$ be a closed Riemannian surface and $(M,h)$ a Riemannian manifold. The \emph{energy} of a smooth map $f: \Sigma \to M$ is
\[\E(f) = \frac{1}{2}\int_\Sigma \Vert \d f\Vert^2 \dvol_g~,\]
where the norm $\Vert \d f\Vert$ is computed using the tensor product metric on $\T^* \Sigma\otimes f^*\T M$. The map $f$ is called \emph{harmonic} if it is a critical point of the energy, namely if for any smooth variation $(f_t)_{t\in (-\epsilon,\epsilon)}$ we have
\[\left. \frac{\d}{\d t}\right\vert_{t=0} \E(f_t)=0~.\]

Since $\Sigma$ has dimension $2$, the energy $\E(f)$ only depends on the conformal class of $g$, and so harmonic maps can be defined on the underlying Riemann surface $X=(\Sigma,[g])$: a smooth map $f: X \to (M,h)$ is harmonic if and only if the $\mathbb C$-linear part $\partial f$ of its differential is holomorphic.

The notion of harmonic maps is intimately linked with the theory of branched minimal immersions. The following result due to Robert Gulliver, Robert Osserman and Halsey Royden \cite{Gulliver:1973} is  an extension of a classical observation of James Eells and Joseph Sampson \cite{EellsSampson} -- see also Jonathan Sacks and Karen Uhlenbeck \cite[Theorem 1.6]{Sacks:1982aa}.

\begin{proposition}\label{prop:EellsSampson}
A map $f$ from $X$ to $(M,h)$ is a branched minimal immersion if and only if it is harmonic and conformal.
\end{proposition}

The lack of conformality of a harmonic map $f: X \to (M,h)$ is encoded in its \emph{Hopf differential}, which is defined by
\[\Hopf(f)= h^{\mathbb C}(\partial f,\partial f)~,\]
where $h^{\mathbb C}$ is the $\mathbb C$-linear extension of $h$. Observe that since $\partial f$ is holomorphic, $\Hopf(f)$ is a holomorphic quadratic differential on $X$, that is, an element of $H^0(K_X^2)$. The following proposition is a direct consequence of Proposition \ref{prop:EellsSampson}:

\begin{proposition}\label{prop:ZeroHopf}
A harmonic map $f$ from $X$ to $(M,h)$ is a branched minimal immersion if and only if $\Hopf(f)=0$.
\end{proposition} 

%Let us observe that the proof of this result is local, so we do not need to assume that $X$ is closed.

\subsection{The equivariant case} In this paper, every group action on a Riemannian manifold will be by isometries.

\begin{definition}\label{def:Good}
Let $\Gamma$ be a discrete group acting cocompactly on the hyperbolic plane $\hh$ and $\rho$ be a representation of $\Gamma$ into $\mathsf{U}(n)$. The representation $\rho$ has \emph{finite energy} if there exists a $\rho$-equivariant map from $\hh$ to the sphere $\bfS^{2n-1}$.
\end{definition}

Observe that for $n$ greater than 1 and $\Gamma$ torsion free, any unitary representation has  finite energy. However, if $\Gamma$ has torsion, a necessary condition for $\rho$ to be of finite energy is that any torsion element has a fixed point in the sphere.

Given a smooth equivariant map $f$, the pointwise norm $\Vert \d f\Vert$ is $\Gamma$-invariant, so it descends to a function on $\Gamma\backslash \hh$. Thus, if $\Gamma$ is cocompact, possibly with torsion, one can define the energy $\E(f)$ by integrating over a fundamental domain $D$ for the action of $\Gamma$:
\[\E(f) = \frac{1}{2}\int_D\Vert\d f\Vert^2\dvol_g ~.\]
Subsequently, let
\[\E(\rho) \defeq \inf \left\{\E(f)~\vert~ f: \hh\to \bfS^{2n-1}\text{ is smooth and $\rho$-equivariant}\right\}~.\]
%Note that we do not care whether $\Gamma$ has torsion or not.
A smooth equivariant map $f$ will be called \emph{energy minimizing} if $\E(f)=\E(\rho)$.

In \cite[Theorem 1.7]{song}, Song obtained  an extension of a classical result of Sacks and Uhlenbeck \cite{SacksUhlenbeck}, in the spirit of \cite{Corlette:1988,Donaldson:1987,Labourie:1991}. We actually need a  version of this result when torsion could be present:

\begin{theorem}[\sc Sacks--Uhlenbeck, Song]\label{thm:SacksUhlenbeck2}
Let $\Gamma$ be a discrete group acting cocompactly on the hyperbolic plane $\hh$. Let  $\rho$ be a finite energy unitary representation of $\Gamma$ in $\U(n)$. Then, there exists a $\rho$-equivariant energy minimizing harmonic map $\psi$ from $\hh$ to $\bfS^{2n-1}$.
\end{theorem}

\begin{proof}\label{app:karen}
Consider the space
\[\mathcal E_{\rho} = \left\{ f: \hh \to \bfS^{2n-1}~\vert ~f \text{ is $C^\infty$ and $\rho$-equivariant} \right\}~.\]
By assumption $\mathcal E_{\rho}$ is nonempty so we can find a sequence $(f_j)_{j\in \N}$ such that $$\lim \E(f_j) = \inf \{\E(f)~\vert~f\in \mathcal E_{\rho}\}~.$$

Let $D$ be the discrete set of points in $\hh$ whose stabilizer in $\Gamma$ is nontrivial. Let $B$ be a ball in $\hh$ whose iterates under $\Gamma$ cover $\hh$. Embed $\mathcal E_{\rho}$ in the Sobolev space $W^{1,2}(B,\mathbb R^{2n})$ (recall that the energy is defined on $W^{1,2}(B,\mathbb R^{2n})$). The sequence $(f_j)_{j\in \N}$ is now bounded in $W^{1,2}(B,\mathbb R^{2n})$ hence weakly converges by Banach--Alaoglu to $\psi$.  By Rellich--Kondrachov Theorem, the embedding of $W^{1,2}(B,\mathbb R^{2n})$ into $L^2(B,\mathbb R^{2n})$ is compact, and so $(f_j)_{j\in \N}$ also strongly converges in $L^2(B,\mathbb R^{2n})$ to $\psi$.  From that we observe that $\psi$ is almost everywhere $\rho$-equivariant and takes value in $\bfS^{2n-1}$, hence can be equivariantly extended to a map defined on $\hh$.

By lower semi-continuity of the energy, one obtains that $\E(\psi)\leq \lim \E(f_j)$, and so $\psi$ minimizes locally the energy, away from $D$. By \cite[Section 8.4.3]{Evans:2010aa} one gets that $\psi$ is a weak solution to the harmonic equation on $\hh\setminus D$. Finally, a classical result of Frédéric Hélein \cite{Helein:1991aa} implies that $\psi$ is a strong solution on $\hh\setminus D$. 

To extend along the discrete set $D$, observe that any point $p$ in $D$ has finite stabilizer in $\Gamma$. In particular, $\psi$ has finite energy on a neighborhood $U$ of $p$. By a result of Sacks--Uhlenbeck \cite[Theorem 3.6]{SacksUhlenbeck}, $\psi$ extends to a smooth harmonic map on $U$.
\end{proof}

\subsection{Strong convergence}

Given a discrete group $\Gamma$, the \emph{left regular representation} of $\Gamma$ is the unitary representation
\[\lambda_\Gamma: \Gamma \longrightarrow \U\big(\ell^2(\Gamma,\mathbb C)\big)\]
defined by $\left(\lambda_\Gamma(\gamma)(f)\right)(x)=f(\gamma^{-1}x)$.

A sequence $(\rho_j)_{j\in \mathbb N}$ of representations with $\rho_j: \Gamma \to \U(N_j)$ \emph{strongly converges} to a  representation  $\rho$ if
\[\forall f\in \mathbb C[\Gamma]~,\lim_{j\to\infty}\Vert \rho_j(f)\Vert = \Vert \rho(f)\Vert~,\]
where the norms are the operator norms. Recall that the operator norm of a bounded linear operator $A$ is 
\[ \Vert A \Vert = \sup\{\Vert A(u)\Vert ~,~ \Vert u\Vert=1\}~.\]
We say that two representations $\rho_1$ and $\rho_2$ are \emph{weakly equivalent} if 
\[\forall f\in \mathbb C[\Gamma]~,\Vert \rho_1(f)\Vert = \Vert \rho_2(f)\Vert~.\]
Note that if a sequence $(\rho_j)_{j\in \mathbb N}$ of representations with $\rho_j: \Gamma \to \U(N_j)$ \emph{strongly converges} to  a  representation  $\lambda_1$, and  furthermore $\lambda_1$ and  $\lambda_2$ are weakly equivalent, then $(\rho_j)_{j\in \mathbb N}$ of also strongly converges to $\lambda_2$.

We finally say  a sequence $(\rho_j)_{j\in \mathbb N}$ {\em virtually strongly converges with respect to  a subgroup $\Gamma_0$} to the regular representation if there is a finite index subgroup $\Gamma_0$ such that the sequence restricted to $\Gamma_0$ strongly converges to $\lambda_{\Gamma_0}$

The following is provided by Lars Louder and Michael Magee \cite{LMH}:

\begin{theorem}[\sc Louder--Magee]\label{thm:StrongConv}
Let $\Gamma$ be the fundamental group of a closed connected oriented surface of genus at least $2$. There exists a sequence $(\rho_j)_{j\in \mathbb N}$ of unitary representations of $\Gamma$ in $\U(n_j)$ of finite image, that strongly converges to the regular representation.
\end{theorem}

\subsection{A result of Antoine Song}

Recall the following theorem of Song \cite[Theorem 0.4]{song}.

\begin{theorem}[\sc Song's Convergence Theorem]\label{thm:Song}
Let $X$ be a closed Riemann surface. Let $(\rho_j)_{j\in \mathbb N}$ be a sequence, where   $\rho_j$ is a representation of  $\pi_1(X) $ in $\U(N_j)$. If $(\rho_j)_{j\in \mathbb N}$ strongly converges to the regular representation  
 then,
 $$\lim_{j\to\infty}\E(\rho_j)=\tfrac{\pi}{4}\vert \chi(X)\vert\ .
 $$
 Moreover if  $\psi_j$ from $\hh$ to $ \bfS^{2N_j-1}$ is a $\rho_j$-equivariant energy minimizing harmonic map, then
\[\lim_{j\to\infty} \psi_j^*\g_{\bfS^{2N_j-1}}=\tfrac{1}{8}\g_{\hh}~,\]
where $\g_{\bfS^k}$ is the round metric on the sphere $\bfS^k$, $\g_{\hh}$ is the hyperbolic metric on $\hh$ and the convergence is in the $C^\infty$ topology.
\end{theorem}

We note that the convergence of the energy in the above theorem can be thought of as a non-linear version of \cite{HideMagee}.

The following rigidity result \cite[Corollary 2.4]{song}, strenghtened in \cite{Caniato:2025aa}, is crucial.

\begin{theorem}[\sc Song's rigidity theorem]\label{thm:SongRigidity}
Let $X$ be a closed Riemann surface, $\cH$ be a Hilbert space and $\rho$ be a representation of $\pi_1(X)$ into $\mathsf U(\cH)$ which is weakly equivalent to the regular representation of $\pi_1(X)$. Any $\rho$-equivariant map $\varphi$ from $\hh$ to the unit sphere in $\cH$ satisfies
\[\E(\varphi)\geq \tfrac{\pi}{4}\vert \chi(X)\vert\]
and equality holds if and only if
\[\varphi^* \g_{\bfS(\cH)} = \tfrac{1}{8}\g_\hh~,\]
where $\g_{\bfS(\cH)}$ is the round metric on the unit sphere of $\cH$.
\end{theorem}

Given $\Gamma$ a discrete group (possibly with torsion) acting cocompactly on $\hh$, one obtains an \emph{orbifold Riemann surface} $X=\Gamma\backslash\hh$. By Selberg  Lemma, $\Gamma$ admits a torsion free normal subgroup $\Gamma_0$ of finite index $[\Gamma:\Gamma_0]$. In particular, the Riemann surface $X_0=\Gamma_0\backslash\hh$ is smooth and is a branched cover of $X$ of order $[\Gamma:\Gamma_0]$. The \emph{orbifold Euler characteristic} $\chi_{\scriptscriptstyle{\mathrm{orb}}}(X)$ of $X$ is defined by 
\[\chi_{\scriptscriptstyle{\mathrm{orb}}}(X):= \frac{1}{[\Gamma:\Gamma_0]}\chi(X_0)~.\]

We now prove the following consequence of Song's result in the torsion case:

\begin{theorem}\label{thm:SongTorsion}
Let $\Gamma$ be a discrete group acting cocompactly on $\hh$. Let $(\rho_j)_{j\in \mathbb N}$ be a sequence of representations, of $\Gamma$ into $\U(N_j)$, assume that 
\begin{align}
\lim_{j\to \infty} \E(\rho_j) = \tfrac{\pi}{4}\vert \chi_{\scriptscriptstyle{\mathrm{orb}}}(\Gamma\backslash \hh)\vert ~.\label{eq:limrho-torsion}
\end{align}

Assume furthermore that $(\rho_j)_{j\in\mathbb N}$  virtually strongly converges with respect to a torsion free subgroup to the regular representation.
Let $u_j$ be a $\rho_j$-equivariant  harmonic map from $\hh$ to $\bfS^{2N_j-1}$. Assume that
$$\lim_{j\to \infty} \E(u_j) = \tfrac{\pi}{4}\vert \chi_{\scriptscriptstyle{\mathrm{orb}}}(\Gamma\backslash \hh)\vert\ .$$ 
	Then
\[\lim_{j\to \infty} u_j^* \g_{\bfS^{2N_j-1}} = \tfrac{1}{8}\g_\hh\]
where the convergence is in the $C^\infty$ topology.
\end{theorem}

\begin{proof}
Let $\Gamma_0$ be the torsion free subgroup of $\Gamma$ such that $(\rho^0_j)_{j\in\mathbb N}$  strongly converges to the regular representation, where $\rho^0_j$ is the restriction of $\rho_j$ to $\Gamma_0$.

For any $j$ in $\mathbb N$, let $u^0_j$ be the harmonic map $u_j$ seen as a $\rho^0_j$-equivariant map. Observe that
\[\E(u^0_j) = \frac{1}{2}\int_{\Gamma_0\backslash \hh} \Vert \d u^0_j\Vert^2 \dvol_{\hh} = \frac{[\Gamma:\Gamma_0]}{2}\int_{\Gamma\backslash \hh} \Vert \d u_j\Vert^2\dvol_{\hh} = [\Gamma:\Gamma_0] \ \E(u_j)~.\]
In particular, $u^0_j$ might not be energy minimizing, but we still have $$
\lim_{j\to\infty} \E(u^0_j) = [\Gamma:\Gamma_0]\left(\frac{\pi}{4}\vert \chi_{\scriptscriptstyle{\mathrm{orb}}}(\Gamma\backslash \hh)\vert \right)= \frac{\pi}{4}\vert \chi(\Gamma_0\backslash \hh)\vert\ .$$

We now follow the proof of \cite[Theorem 0.4]{song} to apply the conclusion of Theorem \ref{thm:Song} to the closed Riemann surface $\Gamma_0\backslash \hh$. By embedding each sphere totally geodesically into a fixed infinite dimensional sphere $\bfS(\cH)$ in a Hilbert space $\cH$, we can consider the sequence $(u^0_j)_{j\in\mathbb N}$ as taking value in a fixed Riemannian manifold of infinite dimension. For any $j$, one can find a unitary transformation $h_j$ in $\mathsf U(N_j)$ such that the 1-jet of $h_j \circ u^0_j$ at a fixed point $x$ in  $ \hh$ remains in a fixed finite dimensional compact submanifold.

By the upper bound on the energy of $(u^0_j)_{j\in \N}$, one can apply Sacks--Uhlenbeck compactness result \cite[Theorem 4.4]{SacksUhlenbeck} (which follows from $\epsilon$-regularity, see \cite[Theorem 1.8]{song}): there is a discrete subset $D$ of $\hh$ such that $(u^0_j)_{j\in \N}$ converges $C^\infty$ on any compact of $\hh\setminus D$ to a harmonic map $u^0_\infty$ of locally finite energy defined on $\hh$. By \cite[Theorem 4.6]{SacksUhlenbeck}\footnote{The original proof of Sacks--Uhlenbeck holds for compact finite dimensional range. Nevertheless, its proof only relies on the $\epsilon$-regularity which holds uniformly for all spheres (as pointed out in \cite[Section 1.2]{song}).}, at any point of $D$ where $C^\infty$-convergence does not hold, one obtains a \emph{bubble}, that is, a nonconstant harmonic map $\beta$ from $\bfS^2$ to $\bfS(\cH)$ such that 
\begin{equation}\label{eq:bubble}
\E(u^0_\infty) + \E(\beta) \leq \lim_{j\to\infty} \E(u^0_j)=\frac{\pi}{4}\vert\chi (\Gamma_0\backslash \hh)\vert~.
\end{equation}
By \cite[Proof of Proposition 4.4]{song}, the map $u^0_\infty$ is equivariant under a representation weakly equivalent to the regular representation of $\Gamma_0$. Since $\beta$ has positive energy, equation \eqref{eq:bubble} contradicts the first statement of Theorem \ref{thm:SongRigidity}. It follows that $(u^0_j)_{j\in \N}$ converges $C^\infty$ on any compact of $\hh$ to $u^0_\infty$ and $\E(u^0_\infty)=\frac{\pi}{4}\vert\chi (\Gamma_0\backslash \hh)\vert$. The second statement of Theorem \ref{thm:SongRigidity} gives the result.
\end{proof}

\begin{remark}
We will see in Theorem  \ref{theo:InducedRep2} that if the sequence $(\rho_j)_{j\in \mathbb N}$ above is obtained as the sequence of induced representations of a strongly convergence sequence, then assumption \eqref{eq:limrho-torsion} is automatically satisfied.
\end{remark}

\section{Induced representation}

Let $\Gamma$ be a discrete group acting cocomplactly on $\hh$ and $\Gamma_0$ a finite index normal subgroup. In this section, we apply the construction of the \emph{induced representation} to our setting. We recall the basics of the theory in Section \ref{sec induced}, more details  can found in the book of Robert Zimmer \cite{Zimmer:1984aa}.  

Given a Hilbert space $\cH_0$ and a representation $\rho_0$ of $\Gamma_0$ into $\mathsf{U}(\cH_0)$, the induced representation produces a new Hilbert space $\cH_0^\ind$ together with a representation $\rho_0^\ind$ of $\Gamma$ into $\mathsf{U}(\cH_0^\ind)$ satisfying the following

\begin{theorem}[\sc Properties of the induced representation]\label{thm:InducedRep}
Let $\Gamma$ be a discrete group acting cocompactly on $\hh$, and $\Gamma_0$ be a finite index normal subgroup of $\Gamma$. Given a representation $\rho_0$ of $\Gamma_0$ into $\mathsf U(\cH_0)$ for some Hilbert space $\cH_0$, the induced representation $\rho_0^\ind$ of $\Gamma$ into $\mathsf U(\cH_0^\ind)$ satisfies the following properties
\begin{enumerate}
	\item  \label{it:1-induced} if $\cH_0$ has finite dimension, then $\cH_0^\ind$ has finite dimension.
	\item  \label{it:2-induced} If $\rho_0$ has finite image, then $\rho_0^\ind$ has finite image.
	\item \label{it:3-induced} If $\rho_0$ has finite energy in the sense of Definition \ref{def:Good}, then $\rho_0^\ind$ has finite energy and 
	\[\E(\rho_0^\ind)\leq \frac{1}{[\Gamma:\Gamma_0]}\E(\rho_0)~.\]
	\item  \label{it:4-induced} If $(\rho_j)_{j\in \mathbb N}$ is a sequence of unitary representations of $\Gamma_0$ that strongly converges to the regular representation, then the sequence $(\rho_j^\ind)_{j\in\mathbb N}$ of induced representations of $\Gamma$ virtually strongly converges with respect to $\Gamma_0$ to the regular representation.
	%\footnote{Actually, it is well known that in this situation  $(\rho_j^\ind)_{j\in\mathbb N}$ strongly converges to the regular representation. If the sequence strongly converges, then it virtually strongly converges with respect to any finite index subgroup. The direct proof presented here is for the sake of the reader.}
\end{enumerate}
\end{theorem}
As we will prove in Sections \ref{sec induced} and \ref{sec induced 2}, the previous assertion follows from  elementary considerations.

\begin{theorem}[\sc Song's Theorem for induced representations]\label{theo:InducedRep2}

Let $\Gamma$ be a discrete group acting cocompactly on $\hh$ and $\Gamma_0$ be a torsion free normal subgroup of finite index. Let $(\rho^0_j)_{j\in \mathbb N}$ be a sequence of unitary representations of $\Gamma_0$ in  $\U(n_j)$ strongly converging to the regular representation.  
Let $(\rho^{\ind}_j)_{j\in \mathbb N}$ be the sequence of unitary representations of $\Gamma$ in  $\U(N_j)$ induced by $(\rho^0_j)_{j\in \mathbb N}$.  
Then 
$$
\lim_{j\to\infty}{\E(\rho^{\ind}_j)}=\frac{\pi}{4}\left|\chi_{\scriptscriptstyle{\mathrm{orb}}}(\Gamma\backslash\hh)\right|~.
$$	
Moreover, if we choose for each $j$ a $\rho^{\ind}_j$-energy minimizing map  $u_j$ then
\[\lim_{j\to \infty} u_j^* \g_{\bfS^{2N_j-1}} = \frac{1}{8}\g_\hh\]
where the convergence is in the $C^\infty$ topology.
\end{theorem}
We postpone the proof of Theorem \ref{theo:InducedRep2} until Section \ref{proof of induced}.

\subsection{The induced representation}\label{sec induced}
We now describe more precisely the induced representation. Let $\Gamma_0$ be a finite index normal subgroup of $\Gamma$. Assume that we have a representation $\rho_0$ from $\Gamma_0$ to $\U(\cH_0)$ for $\cH_0$ some Hilbert space. Define the vector space
$$
\cH^\ind_0:=\big\{f:\Gamma\to \cH_0\mid \forall \gamma_0 \in\Gamma_0,~x\in \Gamma \ \  f(\gamma_0x)=\rho_0(\gamma_0)f(x)\big\}\ .
$$
Note that $\cH^\ind_0$ is a Hilbert space when equipped with the scalar product
\[
\langle f,g\rangle_{\cH^\ind_0} =\frac{1}{[\Gamma:\Gamma_0]}\sum_{\eta\in \Gamma/\Gamma_0} \langle f(\eta),g(\eta)\rangle_{\cH_0}\ .
\]
The {\em induced representation} $\rho^\ind_0$ from $\rho_0$ is the unitary representation of $\Gamma$ into $\mathsf{U}(\cH^\ind_0)$ given    by
\[
\rho_0^\ind(\gamma)f\ :\  \mapping{\Gamma}{\cH_0\ ,} 
{\eta}{f(\eta\gamma)
\ . }\]

Item \eqref{it:2-induced} of Theorem \ref{thm:InducedRep} follows from

\begin{lemma}
If $\rho_0$ has finite image then so has $\rho^\ind_0$.
\end{lemma}
\begin{proof}
	It is enough to prove the result when $\rho_0$ is trivial. In that case, $\cH^\ind_0$ is the space of maps from $\Gamma$ to $\cH_0$ which are left $\Gamma_0$-invariant, that is, the space of $\cH_0$-valued maps on $S=\Gamma/\Gamma_0$. Since $\Gamma_0$ acts trivially on $S$, for all $\gamma_0$ in $\Gamma_0$, we have $\rho^\ind_0(\gamma_0\gamma)=\rho^\ind_0(\gamma)$. Thus $\Gamma_0$ is in the kernel of $\rho^\ind_0$. Hence $\rho^\ind_0$ have finite image.
	 \end{proof}
	 
The proofs of items \eqref{it:1-induced} and \eqref{it:4-induced} follow directly from

\begin{lemma}\label{lem:ind-ident}
For $s$ in $\Gamma$ and  $h$ in $\cH_0$, let $\rho_s(\gamma_0)h=\rho_0(s\gamma_0 s^{-1})h$. Let us fix a set $S$ in $\Gamma$, such that $S$ represents $\Gamma/\Gamma_0$. Let $\hat\rho_0^\ind$ be the restriction of $\rho_0^\ind$ to $\Gamma_0$. Then $\hat\rho_0^\ind$ is isomorphic, as $\Gamma_0$-representations, to $\bigoplus_{s\in S} \rho_s$. In particular, for any element $f$ in $\mathbb C[\Gamma_0]$ we have
\[\Vert \hat\rho_0^\ind(f)\Vert = \Vert \rho_0(f)\Vert~.\]
where the norm is the operator norm.
\end{lemma}

\begin{proof}
Define 
\[
\theta: \mapping{\cH_0^\ind }{\bigoplus_{s\in S} \cH_0\ , }{f}{\left(f(s)\right)_{s\in S}\ .}\]

By construction, $\theta$ is an isometry, where the Hilbert product on the right-hand side is the direct sum scalar product (in particular, the diagonal embedding of $\cH_0$ is an isometry). Moreover, we have

\begin{align*}
\theta\left(\hat\rho_0^\ind(\gamma_0) f\right) & = \left(f(s\gamma_0)\right)_{s\in S} \\
& = \left( f\left((s\gamma_0 s^{-1})s\right) \right)_{s\in S}\\
& = \left(\rho_s(\gamma_0)f(s)\right)_{s\in S}
\end{align*}

In particular, $\theta$ intertwines $\hat\rho_0^\ind$ and $\bigoplus_{s\in S} \rho_s$.

The last statement follows from the fact that, for any endomorphism $\varphi$ of a Hilbert space $\cH$, the operator norm of $\varphi \oplus \varphi$ on $\cH \oplus \cH$, equipped with the direct sum scalar product, satisfies $\Vert \varphi \oplus \varphi\Vert = \Vert\varphi\Vert$.
\end{proof}

\subsection{Induced equivariant map}\label{sec induced 2}

We now prove item \eqref{it:3-induced} of Theorem \ref{thm:InducedRep}, namely
\begin{proposition}\label{prop item 3}
	Let $\Gamma$ be a discrete group acting cocompactly on $\hh$ and $\Gamma_0$  be a finite index normal subgroup. If $\rho_0$ is a unitary representation of $\Gamma_0$ of finite energy (in the sense of Definition \ref{def:Good}), then  $\rho^\ind_0$ has finite energy  and
\[\E(\rho^\ind_0)\leq \frac{\E(\rho_0)}{[\Gamma:\Gamma_0]}  \ .\]
\end{proposition}

Consider a $\rho_0$-equivariant map $\varphi_0$ from $\hh$ to $\bfS(\cH_0)$. Thus for any $h$ in $\hh$, $\varphi_0(h)$ is an element of $\cH_0$ of norm $1$. The \emph{induced equivariant map} is the map $\varphi_0^\ind$ from $\hh$ to $\cH_0^\ind$, where $\varphi_0^\ind(h)$ is defined by
$$
\mapping{\Gamma}{\cH_0}{\eta}{\varphi_0(\eta h)\ .}
$$
Observe that since $\varphi_0(h)$ has norm $1$, then
$$
\Vert \varphi_0^\ind(h)\Vert^2_{\cH^\ind_0}=\frac{1}{[\Gamma:\Gamma_0]}\sum_{\eta\in \Gamma/\Gamma_0} \Vert \varphi_0(\eta h)\Vert_{\cH_0}^2=1\ ,
$$
and so $\varphi_0^\ind$ takes value in $\bfS(\cH_0^\ind)$. The proof of Proposition \ref{prop item 3} is a direct consequence of the following

\begin{lemma}
The map $\varphi_0^\ind$ is $\rho_0^\ind$-equivariant and its energy satisfies 
\[\E(\varphi_0^\ind)= \frac{\E(\varphi_0)}{[\Gamma:\Gamma_0]}~.\]
\end{lemma}

\begin{proof}
For $\gamma$ in $\Gamma$, the element $\varphi_0^\ind(\gamma h)$ is given by
$$
\mapping{\Gamma}{\cH_0}{\eta}{\varphi_0\left(\eta(\gamma h)\right)\ .}
$$
On the other hand $\rho_0^\ind(\gamma)(\varphi_0^\ind(h))$ is given 
$$
\mapping{\Gamma}{\cH_0}{\eta}{\varphi_0((\eta \gamma) h)\ .}
$$
Thus $\rho_0^\ind(\gamma)\varphi_0^\ind(h)=\varphi_0^\ind(\gamma h)$ which means that $\varphi_0^\ind$ is $\rho^\ind_0$-equivariant. Then, for any $h \in \hh$, using the identification of Lemma \ref{lem:ind-ident}, we have
\begin{align*}
\varphi_0^\ind(h)&=\sum_
{\eta\in\Gamma/\Gamma_0}\varphi_0\circ\eta (h)\ ,
\\	\Vert\d_h \varphi_0^\ind \Vert^2 &=\frac{1}{[\Gamma: \Gamma_0]}\sum_{\eta\in \Gamma/\Gamma_0} \Vert\d_{\eta h} \varphi_0\Vert^2\ ,
\end{align*}
where in the last equation we use that $\eta$ is an isometry of $\hh$.  From this it follows that
\begin{align*}
\E(\varphi_0^\ind) & =  \frac{1}{2}\int_{\Gamma\backslash\hh} \Vert \d_h \varphi_0^\ind\Vert^2 \dvol_{\hh}(h) \\
% &=  \frac{1}{2[\Gamma:\Gamma_0]}\int_{\Gamma_0\backslash\hh} \Vert \d_h \varphi_0^\ind\Vert^2 \dvol_{\hh}(h) \\
& =  \frac{1}{2[\Gamma:\Gamma_0]}\int_{\Gamma\backslash\hh} \Big(\sum_{\eta\in \Gamma/\Gamma_0} \Vert \d_{\eta h} \varphi_0\Vert^2 \Big) \dvol_{\hh}(h) \\
& =  \frac{1}{2[\Gamma:\Gamma_0]} \int_{\Gamma_0\backslash\hh} \Vert \d_h \varphi_0\Vert^2 \dvol_{\hh}(h)  \\ & =  \frac{1}{[\Gamma:\Gamma_0]}\E(\varphi_0)~,  \\
\end{align*}
where the third equality follows from the fact that if $D_\Gamma$ is a fundamental domain for the action of $\Gamma$ on $\hh$, its orbit under $ \Gamma/ \Gamma_0$ is a fundamental domain for the action of $\Gamma_0$. The result follows.
\end{proof}
\subsection{Proof of Theorem \ref{theo:InducedRep2}}\label{proof of induced}
Let $\hat\rho_j$ be the restriction of $\rho^{\ind}_j$ to $\Gamma_0$.  By the fourth item of Theorem \ref{thm:InducedRep}, $(\hat\rho_j)$ converges strongly to the regular representation of $\Gamma_0$. In particular the second hypothesis of Theorem \ref{thm:SongTorsion} holds. Then, by applying Theorem \ref{thm:Song} twice, we obtain
\begin{align}
\lim_{j\to\infty}\E(\hat\rho_j)= \frac{\pi}{4}\vert\chi(\Gamma_0\backslash\hh)\vert\ ,	 \crcr
\lim_{j\to\infty}\E(\rho^0_j)= \frac{\pi}{4}\vert\chi(\Gamma_0\backslash\hh)\vert\ .	\label{ineq:ener1} 
\end{align}

By item \eqref{it:3-induced} of  Theorem \ref{thm:InducedRep} the induced representations have finite energy. By Theorem \ref{thm:SacksUhlenbeck2}, for each $j$ there exists a $\rho_j^\ind$-equivariant energy-minimizing map $u_j$. In particular, seeing $u_j$ as a $\hat\rho_j$-equivariant map (denoted by $\hat u_j$), we get
\[\E(\hat u_j) = [\Gamma:\Gamma_0] \E(u_j) = [\Gamma:\Gamma_0]\E(\rho_j^\ind)~.\] 
Since $\hat u_j$ is $\hat\rho_j$-equivariant, we have $\E(\hat \rho_j)\leq \E (\hat u_j)$. Using the above equality, this gives
\[\frac{1}{[\Gamma:\Gamma_0]}\E(\hat \rho_j) \leq \E(\rho_j^\ind) ~.\]
On the other hand, item \eqref{it:3-induced} of Theorem \ref{thm:InducedRep} gives
\[\E(\rho_j^\ind)\leq \frac{1}{[\Gamma:\Gamma_0]} \E(\rho^0_j)~.\]
Taking the limit as $j$ goes to infinity together with equation \eqref{ineq:ener1}, we get $\lim \E(\rho_j^\ind) = \frac{\pi}{4}\vert \chi_{\scriptscriptstyle{\mathrm{orb}}}(\Gamma\backslash \hh)\vert$, thus proving the first assertion. The sequences $(u_j)_{j\in\mathbb N}$ and $(\rho^\ind_j)_{j\in\mathbb N}$ satisfy all the hypothesis of Theorem \ref{thm:SongTorsion}, thus we get the result. \qed

\section{Proof of Main Theorem}

\subsection{A Riemann surface with a large automorphism group}
We prove the existence of an orbifold Riemann surface whose Teichmüller space is a point. Given a closed Riemann surface $X_0$, we denote by $\Aut(X_0)$ its automorphism group, and by $\pi$ the quotient map from $X_0$ to the orbifold $X \defeq \Aut(X_0) \backslash X_0$. We prove

\begin{proposition}\label{prop:ExistenceOfX}
There exists a closed Riemann surface $X_0$ of genus greater than $1$ with no nonzero $\Aut(X_0)$-invariant element in $H^0(K^2_{X_0})$.
\end{proposition}

Before proving Proposition \ref{prop:ExistenceOfX}, let us first prove a lemma.
\begin{lemma}
Let $q_2$ be a meromorphic quadratic differential on $X \defeq \Aut(X_0) \backslash X_0$ and $D$ be the set of orbifold points of $X$. The pullback $\pi^* q_2$ to $X_0$ is holomorphic if and only if $q_2$ has at most simple poles at $D$.
\end{lemma}

\begin{proof} 
The computation is local: around a singular point in $D$ the map $\pi$ is given by $z \mapsto z^p$ for some positive integer $p$. In this coordinate system, we can write $q_2= z^{-n}f(z)\d z^{\otimes 2}$ where $f$ is a non vanishing holomorphic function. In particular
\[\pi^* q_2 = z^{-pn}f(z^p) (\d z^p)^{\otimes 2} = z^{2p-2-pn} p^2f(z^p)\d z^{\otimes 2}~.\]
So $\pi^*q_2$ is holomorphic if and only if $n\leq 1$, and the result follows.
\end{proof}

\begin{proof}[Proof of Proposition \ref{prop:ExistenceOfX}]

Let $X$ is an orbifold structure on $\mathbf{P}^1$ with three conical points $D$ of order $p$, $q$ and $r$, with $\frac{1}{p}+\frac{1}{q}+\frac{1}{r}<1$. Let $\Gamma$ be the orbifold fundamental group of $X$. By the choice of the order of the conical points, $\Gamma$ is a subgroup of $\mathrm{PSL}(2,\mathbb{R})$, so that $X\simeq \Gamma\backslash \mathbf{H}^2$. By Selberg lemma \cite[Theorem 4.8.2]{morris2015}, $\Gamma$ contains a finite index torsion free normal subgroup $\Gamma_0$. The group $\Gamma_0$ then acts freely on $\mathbf{H}^2$, so that the surface $X_0:=\Gamma_0\backslash\mathbf{H}^2$ is a smooth Riemann surface. Such surface is also closed, since $\Gamma_0$ has finite index in $\Gamma.$ We claim that the Riemann  surface $X_0$ has no nonzero $\Aut(X_0)$-invariant element in $H^0(K^2_{X_0})$. Indeed, remark that, by construction, $\Aut(X_0)$ contains $\Gamma/\Gamma_0$, so that 
the quotient $\Aut(X_0)\backslash X_0$ is the orbifold $X$ with singular points $D$. By the previous lemma, any $\Aut(X_0)$-invariant element of $H^0(K^2_{X_0})$ is the pullback via the quotient map $\pi: X_0 \to X$ of an element of $H^0(K^2_X(D))$.
But $X$ is an orbifold structure on $\mathbf{P}^1$ with three singular points $D$, so that the $\Aut(X_0)$-invariant elements of $H^0(K^2_{X_0})$ are the pullbacks of elements of $H^0(K^2_{\mathbf{P}^1}(D))$. Since $K_{\mathbf{P}^1}^2(D)$ has degree $-1$, the space $H^0(K^2_{\mathbf{P}^1}(D))$ is trivial. 
\end{proof}

\subsection{Proof of the Main Theorem} Let $X_0$ be the Riemann surface obtained in Proposition \ref{prop:ExistenceOfX} and $\Gamma_0$ be its fundamental group. Consider $\Gamma$ the semidirect product of $\Gamma_0$ by $\Aut(X_0)$. Hence $\Gamma/\Gamma_0= \Aut(X_0)$, so $\Gamma$ is a discrete group containing $\Gamma_0$ as a normal subgroup. Remark moreover that $\Gamma \backslash \hh = \Aut(X_0) \backslash X_0$.
 
By Theorem \ref{thm:StrongConv}, there exists a sequence $(\rho^0_j)_{j\in\mathbb{N}}$ of unitary representations of $\Gamma_0$ with finite image and strongly converging to the regular representation. Let $\rho_j^{\ind}$ denote the representation of $\Gamma$ induced from $\rho^0_j$.
Applying the Theorem \ref{thm:SacksUhlenbeck2}, we obtain a $\rho_j^{\ind}$-equivariant map $u_j$ that is harmonic and energy minimizing.

Since the Hopf differential of $u_j$ is a $\Gamma$-invariant holomorphic quadratic differential, it vanishes by Proposition \ref{prop:ExistenceOfX}. Proposition \ref{prop:ZeroHopf} thus implies that each $u_j$ is a branched minimal immersion. Furthermore, Theorem \ref{theo:InducedRep2} gives
$$\lim_{j\to\infty}u_j^*\g_{\bfS^{2N_j-1}}=\frac{1}{8}\g_{\hh}$$ 
where the convergence is  $C^\infty$. In particular, for $j$ large enough the induced metric $u_j^*\g_{\bfS^{2N_j-1}}$ is positive definite, hence $u_j$ is free of branch points  (see Remark \ref{rem:SingularCurvature}). Thus, for $j$ large enough, $u_j$ is a negatively curved minimal immersion.

Observe also that $\rho^{\ind}_j$ has finite image (by  item \eqref{it:2-induced} of Theorem \ref{thm:InducedRep}), hence, $\ker \rho^{\ind}_j$ has finite index in $\Gamma$. Let $\Gamma_j$ be a torsion-free finite-index subgroup of $\ker \rho^{\ind}_j$, which therefore acts cocompactly on $\hh$, and set
$$X_j=\Gamma_j\backslash\hh\ .$$ 
Then $X_j$ is a closed Riemann surface and $u_j$ is a minimal immersion of $X_j$ in $\bfS^{2N_j-1}$. By construction, the induced metric $u_j^*\g_{\bfS^{2N_j-1}}$ has negative curvature for $j$ large enough. The theorem is proved. \qed

\bibliographystyle{plain}

\begin{thebibliography}{10}

\bibitem{bryant}
Robert~L Bryant.
\newblock Minimal surfaces of constant curvature in {$S^n$}.
\newblock {\em Transactions of the American Mathematical Society},
  290(1):259--271, 1985.

\bibitem{Caniato:2025aa}
Riccardo Caniato, Xingzhe Li, and Antoine Song.
\newblock Area rigidity for the regular representation of surface groups, 2025.

\bibitem{Corlette:1988}
Kevin~R Corlette.
\newblock {Flat $G$-bundles with canonical metrics}.
\newblock {\em Journal of Differential Geometry}, 28(3):361--382, 1988.

\bibitem{Donaldson:1987}
Simon~K. Donaldson.
\newblock {Twisted harmonic maps and the self-duality equations}.
\newblock {\em Proceedings of the London Mathematical Society. Third Series},
  55(1):127--131, 1987.

\bibitem{EellsSampson}
James Eells and Joseph~H Sampson.
\newblock Harmonic mappings of riemannian manifolds.
\newblock {\em American journal of mathematics}, 86(1):109--160, 1964.

\bibitem{Evans:2010aa}
Lawrence~C. Evans.
\newblock {\em Partial differential equations}, volume~19 of {\em Graduate
  Studies in Mathematics}.
\newblock American Mathematical Society, Providence, RI, second edition, 2010.

\bibitem{Gulliver:1973}
Robert Gulliver, Robert Osserman, and Halsey~L Royden.
\newblock {A theory of branched immersions of surfaces}.
\newblock {\em American Journal of Mathematics}, 95:750--812, 1973.

\bibitem{Helein:1991aa}
Fr\'ed\'eric H\'elein.
\newblock R\'egularit\'e{} des applications faiblement harmoniques entre une
  surface et une vari\'et\'e{} riemannienne.
\newblock {\em C. R. Acad. Sci. Paris S\'er. I Math.}, 312(8):591--596, 1991.

\bibitem{HideMagee}
Will Hide and Michael Magee.
\newblock Near optimal spectral gaps for hyperbolic surfaces.
\newblock {\em Ann. of Math. (2)}, 198(2):791--824, 2023.

\bibitem{Labourie:1991}
Fran{\c c}ois Labourie.
\newblock {Existence d'applications harmoniques tordues \`a valeurs dans les
  vari\'et\'es \`a courbure n\'egative}.
\newblock {\em Proc. Amer. Math. Soc.}, 111(3):877--882, 1991.

\bibitem{Lawson:1970aa}
H.~Blaine Lawson, Jr.
\newblock Complete minimal surfaces in {$S^{3}$}.
\newblock {\em Ann. of Math. (2)}, 92:335--374, 1970.

\bibitem{LMH}
Larsen Louder, Michael Magee, and Will Hide.
\newblock Strongly convergent unitary representations of limit groups.
\newblock {\em Journal of Functional Analysis}, 288(6):110803, 2025.

\bibitem{Micallef:1995aa}
Mario~J. Micallef and Brian White.
\newblock The structure of branch points in minimal surfaces and in
  pseudoholomorphic curves.
\newblock {\em Ann. of Math. (2)}, 141(1):35--85, 1995.

\bibitem{morris2015}
Dave~Witte Morris.
\newblock {\em Introduction to arithmetic groups}, volume~2.
\newblock Deductive Press Lieu de publication inconnu, 2015.

\bibitem{SacksUhlenbeck}
Jonathan Sacks and Karen Uhlenbeck.
\newblock The existence of minimal immersions of 2-spheres.
\newblock {\em Annals of mathematics}, 113(1):1--24, 1981.

\bibitem{Sacks:1982aa}
Jonathan Sacks and Karen Uhlenbeck.
\newblock Minimal immersions of closed {R}iemann surfaces.
\newblock {\em Trans. Amer. Math. Soc.}, 271(2):639--652, 1982.

\bibitem{song}
Antoine Song.
\newblock Random harmonic maps into spheres.
\newblock {\em arXiv:2402.10287}, 2024.

\bibitem{yau}
Shing-Tung Yau.
\newblock {\em Seminar on differential geometry}.
\newblock Number 102. Princeton University Press, 1982.

\bibitem{Zimmer:1984aa}
Robert~J. Zimmer.
\newblock {\em Ergodic theory and semisimple groups}, volume~81 of {\em
  Monographs in Mathematics}.
\newblock Birkh\"auser Verlag, Basel, 1984.

\end{thebibliography}

\end{document}